\documentclass[12pt]{amsart}

 \newtheorem{thm}{Theorem}[section]
 \newtheorem{cor}[thm]{Corollary}
 
 \newtheorem{prop}[thm]{Proposition}
 \theoremstyle{definition}
 \newtheorem{defn}[thm]{Definition}
 \theoremstyle{remark}
 
 \numberwithin{equation}{section}

 \newcommand{\norm}[1]{\left\Vert#1\right\Vert}
 
 \newcommand{\C}{\mathbb{C}}

\begin{document}

\title[$N^{p}$-Spaces]
 {$N^{p}$-Spaces}

\author{ Yun-Su Kim }

\address{Department of Mathematics, Indiana University, Bloomington,
Indiana, U.S.A. }

\email{kimys@indiana.edu}


\dedicatory{}

\commby{Daniel J. Rudolph}


\begin{abstract}
We introduce a new norm, called $N^{p}$-norm $(1\leq{p}<\infty)$
on a space $N^{p}(V,W)$ where $V$ and $W$ are abstract operator
spaces. By proving some fundamental properties of the space
$N^{p}(V,W)$, we also discover that if $W$ is complete, then the
space $N^{p}(V,W)$ is also a Banach space with respect to this
norm for $1\leq{p}<\infty$.

\end{abstract}

\maketitle

\section*{Introduction}
For abstract operator spaces $V$ and $W$, a (bounded) linear map
$\phi:V\rightarrow{W}$ provides another linear map
$\phi_{n}:M_{n}(A)\rightarrow{M_{n}(B)}$ defined by
\begin{center}$\phi_{n}((a_{i,j}))=(\phi(a_{i,j}))$\end{center}
where $n=1,2,\cdot\cdot\cdot$. If a sequence
$\{\norm{\phi_{n}}\}_{n=1}^{\infty}$ belongs to $l^{\infty}$, then
$\phi$ is said to be a \emph{completely bounded} map. W.
Stinespring \cite{S} and W. Arveson \cite{A} introduced operator
space theory related to complete boundedness for a map
$\phi:S\rightarrow{B(K)}$ where $S\subset{B(H)}$ and $H$ and $K$
are Hilbert spaces. It also appeared in the 1980s through the
works of G. Wittstock \cite{W, G} and V. Paulsen \cite{P}.
\vskip.1in In this paper, we consider $l^{p}$-norm
($1\leq{p}<\infty$) for the sequence
$\{\norm{\phi_{n}}\}_{n=1}^{\infty}$. Since
$\norm{\phi_{1}}\leq\norm{\phi_{2}}\leq\norm{\phi_{3}}\leq\cdot\cdot\cdot$,
there is no nonzero map $\phi$ such that
$\{\norm{\phi_{n}}\}_{n=1}^{\infty}$ belongs to $l^{p}$. To put it
another way, we define a new norm
\begin{center}\[\norm{\phi}_{p}=\sum_{n=1}^{\infty}\frac{\norm{\phi_{n}}}{n^{p}}\]\end{center}
and study a space $N^{p}(V,W)$ which is a vector space consisting
of all linear maps $\phi:V\rightarrow{W}$ for which
$\norm{\phi}_{p}<\infty$. Then $\norm{\cdot}_{p}$ is a norm,
called $N^{p}$-norm, on the space $N^{p}(V,W)$
$(1\leq{p}<\infty)$. \vskip.2in

 In section 2, we prove fundamental
properties of the space $N^{p}(V,W)(1\leq{p}<\infty)$, for
example, \vskip.1in

\quad (i) $N^{p}(V,W)\subset{N^{q}(V,W)}$ if
$1\leq{p}\leq{q}<\infty$.

\quad (ii) If $\phi:V\rightarrow{W}$ is completely bounded, then
$\phi\in{N^{p}}(V,W)$ for

\quad\quad all $p>1$,\vskip.1in

 \noindent and we characterize a (bounded) linear
map $\phi:V\rightarrow{W}$ by using the space $N^{p}(V,W)$, that
is, the following statements are equivalent :\vskip.1in

\quad (a) $\phi:V\rightarrow{W}$ is a (bounded) linear map.

\quad (b) $\phi\in{N^{p}}(V,W)$ for any $p>2$.

\vskip.1in The main results of this paper are given when $W$ is
complete as follows :\vskip.1in

\quad (i) If $W$ is complete, then the space $B(V,$$W)$ with
$N^{p}$-norm is a

\quad \quad Banach space for $2<p<\infty$.

\quad (ii) If $W$ is complete, then so is $N^{p}(V,W)$ for
$1\leq{p}<\infty$.\vskip.1in


\section{Preliminaries and Notation}
Let $\mathbb{M}_{n,m}(V)$ denote the linear space of $n\times{m}$
matrices with entries from a linear space $V$ and $B(H)$ is the
space of all bounded operators on a Hilbert space $H$ with the
operator norm. We write $\mathbb{M}_{n}(V)=\mathbb{M}_{n,n}(V)$
and if $V=\mathbb{C}$, we let
$\mathbb{M}_{n,m}=\mathbb{M}_{n,m}(\mathbb{C})$. We will denote a
typical element of $\mathbb{M}_{n}(V)$ by $(v_{i,j})$.

\begin{defn}
A \emph{(concrete) operator space} $V$ on a Hilbert space is a
closed subspace of $B(H)$.
\end{defn}
If $V$ is a concrete operator space, then the inclusion
\begin{center}
$\mathbb{M}_{n}(V)\subset\mathbb{M}_{n}(B(H))=B(H^{n})$
\end{center}
provides a norm
$\norm{\hspace{.025in}\cdot\hspace{.025in}}_{\mathbb{M}_{n}(V)}$
on $\mathbb{M}_{n}(V)$, and $M_{n}(V)$ denotes the corresponding
normed space.

We define a \emph{matrix norm}
$\norm{\hspace{.025in}\cdot\hspace{.025in}}$ on a linear space $W$
to be an assignment of a norm
$\norm{\hspace{.025in}\cdot\hspace{.025in}}_{\mathbb{M}_{n}(W)}$
on the matrix space $\mathbb{M}_{n}(W)$ for each
$n\in{\mathbb{N}}$.
\begin{defn}
An \emph{abstract operator space} is a linear space $W$ together
with a matrix norm $\norm{\hspace{.025in}\cdot\hspace{.025in}}$
for which \vskip0.2cm $\bullet$ $\textbf{M1}$
$\norm{v\oplus{w}}_{\mathbb{M}_{m+n}(W)}=$
max$\{\norm{v}_{\mathbb{M}_{m}(W)},\norm{w}_{\mathbb{M}_{n}(W)}\}$
and

$\bullet$ $\textbf{M2}$
$\norm{\alpha{v}\beta}_{\mathbb{M}_{n}(W)}\leq\norm{\alpha}\norm{v}_{\mathbb{M}_{m}(W)}\norm{\beta}$
\vskip0.2cm for all $v\in{\mathbb{M}_{m}(W)}$,
$w\in{\mathbb{M}_{n}(W)}$ and $\alpha\in{M_{n,m}}$,
$\beta\in{M_{m,n}}$.
\end{defn}
By a \emph{linear map} on an abstract operator space $V$ we mean a
bounded linear map defined on $V$. We let $M_{n}(W)$ denote
$\mathbb{M}_{n}(W)$ with the given norm
$\norm{\hspace{.025in}\cdot\hspace{.025in}}_{\mathbb{M}_{n}(W)}$

Given two abstract operator spaces $V$ and $W$ and a linear map
$\phi:V\rightarrow{W}$, we also obtain a linear map
$\phi_{n}:M_{n}(V)\rightarrow{M_{n}(W)}$ defined by
\begin{equation}\label{4}\phi_{n}((v_{i,j}))=(\phi(v_{i,j})).\end{equation}

Since $\phi$ is a bounded map, each $\phi_{n}$ will also be
bounded.

\begin{defn}\cite{P}\label{5}
If $\texttt{sup}_{n}\norm{\phi_{n}}$ is finite, then $\phi$ is
said to be a \emph{completely bounded} map.
\end{defn}
If $\phi$ is completely bounded, then we set
\begin{center}$\norm{\phi}_{cb}=\texttt{sup}_{n}\norm{\phi_{n}},$\end{center}
and $CB(V,W)$ denotes the space of completely bounded maps from
$V$ to $W$.

\vskip0.1cm In Definition \ref{5}, we can also define a
\emph{completely bounded} map as follows : \vskip0.1cm If a
sequence $\{\norm{\phi_{n}}\}_{n=1}^{\infty}$ belongs to
$l^{\infty}$, then $\phi$ is said to be a \emph{completely
bounded} map. Then naturally we can ask when
$\{\norm{\phi_{n}}\}_{n=1}^{\infty}$ belongs to $l^{p}$ where
$1\leq{p}<\infty$. Since
$\norm{\phi_{1}}\leq\norm{\phi_{2}}\leq\norm{\phi_{3}}\leq\cdot\cdot\cdot$,
there is no nonzero map $\phi$ such that
$\{\norm{\phi_{n}}\}_{n=1}^{\infty}$ belongs to $l^{p}$. However,
in the next section, we will introduce a new space, called a
$N^{p}$-space to solve this problem.

\section{The $N^{p}$-Spaces}
Let $V$ and $W$ be abstract operator spaces. For a linear map
$\phi:V\rightarrow{W}$ and $1\leq{p}<\infty$, we introduce a new
norm $\norm{\phi}_{p}$ and a space $N^{p}(V,W)$ in the following
definition.
\begin{defn}\label{1}
Let $V$ and $W$ be abstract operator spaces. If
$\phi:V\rightarrow{W}$ is a linear map and $1\leq{p}<\infty$, then
define a norm
\begin{equation}\label{15}\norm{\phi}_{p}=\sum_{n=1}^{\infty}\frac{\norm{\phi_{n}}}{n^{p}}\end{equation}
and let the space $N^{p}(V,W)$ be a vector space consisting of all
linear maps $\phi:V\rightarrow{W}$ for which
$\norm{\phi}_{p}<\infty$.
\end{defn}
We can easily see that equation (\ref{15}) defines a norm on the
$N^{p}(V,W)$-spaces, and we call $\norm{\phi}_{p}$ the
$N^{p}$-norm of $\phi$.
\begin{prop}\label{2}Let $V$ and $W$ be abstract operator spaces and
$\phi:V\rightarrow{W}$ be a linear map. Then the following
statements are true. \vskip0.1cm $(i).$ If $\phi\in{N^{p}}(V,W)$
for some $1\leq{p}<\infty$, then $\phi\in{N^{q}}(V,W)$ for any
$q\geq{p}$. Thus,
\begin{equation}\label{6}N^{p}(V,W)\subset{N^{q}(V,W)}\end{equation} if $1\leq{p}\leq{q}<\infty$.
\vskip0.1cm $(ii).$ If $\phi:V\rightarrow{W}$ is completely
bounded, then $\phi\in{N^{p}}(V,W)$ for all $p>1$. \vskip0.1cm
$(iii).$ If $\norm{\phi_{n}}\leq{n^{p-1-\epsilon}}$ for some
$\epsilon>0$ and $n=1, 2, 3,\cdot\cdot\cdot,$ then
$\phi\in{N^{p}}(V,W)$.
\end{prop}
\begin{proof}
$(i)$. Suppose that $\phi\in{N^{p}}(V,W)$ and $1\leq{p}\leq{q}$.
For any $n=1, 2, \cdot\cdot\cdot$,\begin{center}
\[\frac{\norm{\phi_{n}}}{n^{q}}\leq\frac{\norm{\phi_{n}}}{n^{p}} .\]\end{center}
It follows that
\begin{equation}\label{7}\norm{\phi}_{q}\leq{\norm{\phi}_{p}}.\end{equation} Since
$\phi\in{N^{p}}(V,W)$, $\norm{\phi}_{p}<\infty$. Thus, from
inequality (\ref{7}), $\norm{\phi}_{q}<\infty$, that is,
$\phi\in{N^{q}(V,W)}$ which proves the inclusion
(\ref{6}).\vskip0.1in

$(ii)$. If $\phi:V\rightarrow{W}$ is completely bounded and
$\norm{\phi}_{cb}=m$, then
\[\norm{\phi}_{p}=\sum_{n=1}^{\infty}\frac{\norm{\phi_{n}}}{n^{p}}\leq{m}\sum_{n=1}^{\infty}\frac{1}{n^{p}}.\]
Since \[\sum_{n=1}^{\infty}\frac{1}{n^{p}}<\infty\] for any
$p>{1}$, we conclude that $\phi\in{N^{p}}(V,W)$ for any $p>{1}$.
\vskip0.1in

$(iii)$. If $\norm{\phi_{n}}\leq{n^{p-1-\epsilon}}$ for some
$\epsilon>0$ and $n=1, 2, 3,\cdot\cdot\cdot,$ then
\[\norm{\phi}_{p}=\sum_{n=1}^{\infty}\frac{\norm{\phi_{n}}}{n^{p}}\leq\sum_{n=1}^{\infty}
\frac{1}{n^{1+\epsilon}}<{\infty}.\] Thus,
$\phi\in{N^{p}}(V,W)$.\vskip0.1cm
\end{proof}
In fact, $\phi\in{N^{p}(V,W)}$ if and only if the sequence
$\{\frac{\norm{\phi_{n}}^{\frac{1}{p}}}{n}\}_{n=1}^{\infty}$
belongs to $l^{p}$ for $1\leq{p}<\infty$. Thus, we can easily see
that Proposition \ref{2} $(i)$ is true.

\begin{prop}\cite{4}\label{16}
If $V$ is an abstract operator space and
$\varphi:V\rightarrow{M_n}$ is a linear map, then
\begin{equation}
\norm{\varphi_n}=\norm{\varphi}_{cb}.\end{equation}
\end{prop}

The set of linear maps from $V$ to $W$ is denoted $B(V, W)$ with
$B(V,$ $V)$ abbreviated $B(V)$.

\begin{cor}
If $V$ is an abstract operator space and
$\varphi:V\rightarrow{M_n}$ is a linear map, then for $p>1$,
$\varphi\in{N^{p}(V,M_{n})}$ and
\begin{equation}
\norm{\varphi}_{p}\leq\norm{\varphi}_{cb}\sum_{n=1}^{\infty}\frac{1}{n^{p}}.\end{equation}
 Furthermore, $CB(V,M_{n})=N^{p}(V,M_{n})=B(V,M_{n})$ for $p>1$.
\end{cor}
\begin{proof}
By Proposition \ref{16},
\[\norm{\varphi_{1}}\leq\norm{\varphi_{2}}\leq\cdot\cdot\cdot\leq\norm{\varphi_{n}}=\norm{\varphi_{n+1}}=\cdot\cdot\cdot=\norm{\varphi}_{cb}.\]
It follows that
\begin{equation}\label{17}
\norm{\varphi}_{p}=\sum_{n=1}^{\infty}\frac{\norm{\varphi_{n}}}{n^{p}}\leq\sum_{n=1}^{\infty}\frac{\norm{\varphi}_{cb}}{n^{p}}
=\norm{\varphi}_{cb}\sum_{n=1}^{\infty}\frac{1}{n^{p}}.
\end{equation}
By Proposition \ref{16} and (\ref{17}), we have
\begin{center}$CB(V,M_{n})=N^{p}(V,M_{n})=B(V,M_{n})$ for $p>1$.
\end{center}
\end{proof}

In Proposition \ref{16}, if $\varphi:V\rightarrow\C$, then
$\norm{\varphi}_{cb}=\norm{\varphi}$ (\cite{P}). Thus, we have the
next Corollary.

\begin{cor}
If $V$ is an abstract operator space, then for $p>1$, each
(bounded) linear functional $f:V\rightarrow{\C}$ belongs to
$N^{p}(V,{\C})$ and
\begin{center}$\norm{f}_{p}=\norm{f}\sum_{n=1}^{\infty}\frac{1}{n^{p}}.
$\end{center}
\end{cor}

\begin{prop}\cite{P2}\label{13}
Let $V$ and $W$ be abstract operator spaces and
$\phi:V\rightarrow{W}$ be a linear map. Then
\begin{equation}\label{14}
\norm{\phi_{n}}\leq{n}\norm{\varphi}.
\end{equation}
\end{prop}

\begin{prop}\label{8}
Let $V$ and $W$ be abstract operator spaces. Then the following
statements are equivalent :\vskip0.1cm (i) $\phi:V{\rightarrow}W$
is a linear map, that is, $\phi\in{B(V, W)}$.\vskip0.1cm (ii)
$\phi\in{N^{p}}(V,W)$ for any $p>2$.
\end{prop}
\begin{proof}
(i) $\rightarrow$ (ii) : By (\ref{14}),
\begin{equation}\label{9}\frac{\norm{\phi_{n}}}{n^{p}}\leq\frac{\norm{\phi}}{n^{p-1}}\end{equation}
for any $p>2$. From inequality (\ref{9}),
\begin{equation}\label{10}\norm{\phi}_{p}=\sum_{n=1}^{\infty}
\frac{\norm{\phi_{n}}}{n^{p}}\leq\sum_{n=1}^{\infty}\frac{\norm{\phi}}{n^{p-1}}\end{equation}
for any $p>2$. Since $\phi\in{B(V, W)}$, we have
\[\sum_{n=1}^{\infty}\frac{\norm{\phi}}{n^{p-1}}<\infty\] for any
$p>2$. From (\ref{10}), we conclude that $\phi\in{N^{p}}(V,W)$ for
any $p>2$.\vskip0.1cm (ii) $\rightarrow$ (i) : Since
$\phi\in{N^{p}}(V,W)$ for any $p>2$,
\[\norm{\phi}_{p}=\sum_{n=1}^{\infty}
\frac{\norm{\phi_{n}}}{n^{p}}=\norm{\phi_{1}}+\sum_{n=2}^{\infty}
\frac{\norm{\phi_{n}}}{n^{p}}<\infty.\] It follows that
$\norm{\phi_{1}}=\norm{\phi}$ is finite. Thus, $\phi\in{B(V, W)}$.
\end{proof}

\begin{thm}\label{3}
Let $V$ and $W$ be abstract operator spaces and $W$ be complete.
Then the space $B(V,$$W)$ with $N^{p}$-norm is a Banach space for
$2<p<\infty$.
\end{thm}
\begin{proof}
By Proposition \ref{8}, it is sufficient to prove that
$N^{p}(V,W)$ is a complete metric space for $2<p<\infty$.

Fix $p$ such that $2<p<\infty$. Let
$\{\varphi_{(k)}\}_{k=1}^{\infty}$ be a Cauchy sequence in
$N^{p}(V,W)$. By the definition of $N^{p}$-norm,
\begin{equation}\label{12}\norm{\varphi_{(n)}-\varphi_{(m)}}\leq\norm{\varphi_{(n)}-\varphi_{(m)}}_{p}.\end{equation}
Let $\epsilon>0$ be given. Since
$\{\varphi_{(k)}\}_{k=1}^{\infty}$ is a Cauchy sequence in
$N^{p}(V,W)$, there exists a sufficiently large integer
$N(\epsilon)>0$ such that if $n>N(\epsilon)$ and $m>N(\epsilon)$,
then $\norm{\varphi_{(n)}-\varphi_{(m)}}_{p}<\epsilon$. From
(\ref{12}), $\{\varphi_{(k)}\}_{k=1}^{\infty}$ is also a Cauchy
sequence in $B(V,W)$.

The completeness of $B(V,W)$ implies that there is
$\psi\in{B}(V,W)$ such that $\lim_{n\rightarrow\infty}$
$\varphi_{(n)}=\psi$. It follows that there exists a sufficiently
large integer $N_{1}(\epsilon)>0$ such that if
$n>N_{1}(\epsilon)$, then $\norm{\varphi_{(n)}-\psi}<\epsilon$.

Since
$\norm{(\varphi_{(n)}-\psi)_k}\leq{k}\norm{\varphi_{(n)}-\psi}$
for any $k=1,2,...$, if $n>N_{1}(\epsilon)$, then
 \begin{center}
 $\norm{\varphi_{(n)}-\psi}_{p}=\sum_{k=1}^{\infty}\frac{\norm{(\varphi_{(n)}-\psi)_k}}{k^{p}}\leq$\end{center}
\begin{center}$\sum_{k=1}^{\infty}\frac{k\norm{\varphi_{(n)}-\psi}}{k^{p}}=\sum_{k=1}^{\infty}
\frac{\norm{\varphi_{(n)}-\psi}}{k^{p-1}}
<\sum_{k=1}^{\infty}\frac{\epsilon}{k^{p-1}}.$\end{center} By our
assumption $(2<p<\infty)$, we have
\[\sum_{k=1}^{\infty}\frac{1}{k^{p-1}}<\infty.\] Since $\epsilon>0$ is arbitrary, we have
lim$_{n\rightarrow\infty}\norm{\varphi_{(n)}-\psi}_{p}=0$. By
Proposition \ref{8}, $\psi\in{N^{p}(V,W)}$ for $2<p<\infty$. Thus,
$N^{p}(V,W)$ is complete for $2<p<\infty$. Note that by
Proposition \ref{8}, $B(V,W)=N^{p}(V,W)$ for $2<p<\infty$. It
follows that the space $B(V,$$W)$ with $N^{p}$-norm is a Banach
space for $2<p<\infty$.
\end{proof}

From Proposition \ref{2}(i) and Proposition \ref{8}, for every
bounded map $\phi:V\rightarrow{W}$, we can find a real number
$r_{\phi}\geq{1}$ defined by
\begin{center}$r_{\phi}=\texttt{inf}\{p:\phi\in{N^{p}(V,W)}\texttt{ and
}1\leq{p}<\infty\}.$\end{center} The number $r_{\phi}$ is called
the \emph{index} of $\phi$.

\begin{thm}
Let $V$ and $W$ be be abstract operator spaces. If $W$ is
complete, then so is $N^{p}(V,W)$ for $1\leq{p}<\infty$.
\end{thm}
\begin{proof}
Suppose that $W$ is complete. Let
$\{\varphi_{(l)}\}_{l=1}^{\infty}$ be a Cauchy sequence in
$N^{p}(V,W)$ for a fixed $p\in[1,\infty)$ and $\epsilon>0$ be
given. Then there is a natural number $N(\epsilon)$ such that for
all natural numbers $n,m\geq{N(\epsilon)}$, we have
\begin{equation}\label{20}\norm{\varphi_{(n)}-\varphi_{(m)}}_{p}=\sum_{k=1}^{\infty}
\frac{\norm{(\varphi_{(n)}-\varphi_{(m)})_{k}}}
 {k^{p}}<\epsilon.\end{equation}
 By inequality (\ref{12}), $\{\varphi_{(l)}\}_{l=1}^{\infty}$ is
also a Cauchy sequence in $B(V,W)$. Since $W$ is complete, so is
$B(V,W)$. It follows that there is a bounded operator
$\varphi\in{B(V,W)}$ such that
\begin{equation}\label{18}\lim_{l\rightarrow\infty}\norm{\varphi_{(l)}-\varphi}=0.\end{equation}

Let $k\in{\{1,2,3,\cdot\cdot\cdot\}}$ be given. It follows from
\eqref{20} that
\[\norm{(\varphi_{(n)}-\varphi_{(m)})_{k}}\leq{k}^{p}\epsilon\]
for
all natural numbers $n,m\geq{N(\epsilon)}$. Thus, for any
$v=[v_{ij}]\in{M_{k}(V)}$,
\begin{equation}\label{26}\norm{(\varphi_{(n)}-\varphi_{(m)})_{k}(v)}
\leq\norm{(\varphi_{(n)}-\varphi_{(m)})_{k}}\norm{v}\leq
{k}^{p}\epsilon\norm{v}\end{equation} if $n,m\geq{N(\epsilon)}$.
Since $\varphi_{(n)}(v_{i,j})$ converges to $\varphi(v_{i,j})$ in
$W$, \eqref{26} implies that
\begin{equation}\label{23}\norm{(\varphi-\varphi_{(m)})_{k}(v)}\leq{k}^{p}\epsilon\norm{v}\end{equation}
if $m\geq{N(\epsilon)}$. It follows from \eqref{23} that
\begin{equation}\label{22}
\norm{(\varphi-\varphi_{(m)})_{k}}\leq{k}^{p}\epsilon,\end{equation}
if $m\geq{N}(\varepsilon)$. Since $\epsilon$ is arbitrary, we have
\begin{equation}\label{24}
\lim_{m\rightarrow\infty}\norm{(\varphi_{(m)}-\varphi)_{k}}=
\lim_{m\rightarrow\infty}\norm{(\varphi-\varphi_{(m)})_{k}}=0
\end{equation}
for any $k\in{\{1,2,3,\cdot\cdot\cdot\}}$. By triangle inequality,
\[\norm{(\varphi_{(n)}-\varphi_{(m)})_{k}}-\norm{(\varphi_{(n)}-\varphi)_{k}}\leq\norm{(\varphi_{(m)}-\varphi)_{k}}\]
and
\[-\norm{(\varphi_{(m)}-\varphi)_{k}}\leq\norm{(\varphi_{(n)}-\varphi_{(m)})_{k}}-\norm{(\varphi_{(n)}-\varphi)_{k}},\]
that is,
\[|\norm{(\varphi_{(n)}-\varphi_{(m)})_{k}}-\norm{(\varphi_{(n)}-\varphi)_{k}}|\leq\norm{(\varphi_{(m)}-\varphi)_{k}},\]
from \eqref{24} we have
\begin{equation}\label{25}
\lim_{m\rightarrow\infty}\norm{(\varphi_{(n)}-\varphi_{(m)})_{k}}=\norm{(\varphi_{(n)}-\varphi)_{k}}
\end{equation}
for any $n$ and $k$ in $\{1,2,3,\cdot\cdot\cdot\}.$

Let $n\geq{N(\epsilon)}$ be given and $\{u_{k}\}_{k=1}^{\infty}$
be a sequence of functions defined on $\{1,2,3,\cdot\cdot\cdot\}$
by
\[u_{k}(m)=\frac{\norm{(\varphi_{(n)}-\varphi_{(m)})_{k}}}{k^{p}}.\]
Then $\sum_{k=1}^{\infty}{u_{k}}$ is uniformly convergent on
$\{1,2,3,\cdot\cdot\cdot\}$, since
$\{\varphi_{(l)}\}_{l=1}^{\infty}$ is a Cauchy sequence in
$N^{p}(V,W).$

Thus, equation \eqref{25} implies that if $n\geq{N(\epsilon)}$,
then
\[\lim_{m\rightarrow\infty}\norm{\varphi_{(n)}-\varphi_{(m)}}_{p}=\lim_{m\rightarrow\infty}
\sum_{k=1}^{\infty}\frac{\norm{(\varphi_{(n)}-\varphi_{(m)})_{k}}}{k^{p}}=\lim_{m\rightarrow\infty}
\sum_{k=1}^{\infty}u_{k}(m)\]\[=\sum_{k=1}^{\infty}\lim_{m\rightarrow\infty}u_{k}(m)=
\sum_{k=1}^{\infty}\lim_{m\rightarrow\infty}\frac{\norm{(\varphi_{(n)}-\varphi_{(m)})_{k}}}{k^{p}}
=\sum_{k=1}^{\infty}\frac{\norm{(\varphi_{(n)}-\varphi)_{k}}}{k^{p}},\]
that is, if $n\geq{N(\epsilon)}$ and $p\in[1,\infty)$,
\begin{equation}\label{29}\lim_{m\rightarrow\infty}\norm{\varphi_{(n)}-\varphi_{(m)}}_{p}
=\norm{\varphi_{(n)}-\varphi}_{p}.\end{equation}

From \eqref{20} and \eqref{29}, we can conclude that
\[\lim_{n\rightarrow\infty}\norm{\varphi_{(n)}-\varphi}_{p}=0,\]
and so $\varphi_{(n)}\rightarrow\varphi$ in $N^{p}$-norm.

Thus, there is a natural number $n_{0}$ such that
\begin{equation}\label{30}\norm{\varphi_{(n_{0})}-\varphi}_{p}\leq\epsilon,\end{equation}
and so by triangle inequality and an inequality (\ref{30}), we
have
\[\norm{\varphi}_{p}=\sum_{k=1}^{\infty}
\frac{\norm{\varphi_{k}}}{k^{p}}\leq\sum_{k=1}^{\infty}\frac{\norm{(\varphi_{(n_0)})_{k}-\varphi_{k}}
+\norm{(\varphi_{(n_0)})_{k}}} {k^{p}}\]
\[=\norm{\varphi_{(n_{0})}-\varphi}_{p}+\norm{\varphi_{(n_{0})}}_{p}\leq\epsilon+\norm{\varphi_{(n_{0})}}_{p}\]
Since $\varphi_{(n_0)}\in{N^{p}(V,W)},$ i.e.,
$\norm{\varphi_{(n_{0})}}_{p}<\infty$,  we have
\[\norm{\varphi}_{p}<\infty.\] Thus, $\varphi\in{N^{p}(V,W)}.$
Therefore, $N^{p}(V,W)$ is complete for $1\leq{p}<\infty$.

\end{proof}

------------------------------------------------------------------------

\bibliographystyle{amsplain}
\bibliography{xbib}
\end{document}